\documentclass[preprint,12pt]{elsarticle}

\usepackage{amsmath}
\usepackage{pgfplots}
\usepackage{amscd}
\usepackage{amssymb}
\usepackage{amsthm}
\usepackage{amsfonts}
\usepackage{stackrel}
\usepackage{tikz-cd}

\def\natu           {\mathbb N}
\def\inte 		{\mathbb Z}
\def\real		{\mathbb R}

\def\comp		{\mathbb C}

\def\R		{\cal R}
\def\L		{\cal L}
\def\F		{\cal F}

\def\lra		{\longrightarrow}

\journal{Journal of Algebra}

\begin{document}

\begin{frontmatter}

\newtheorem{theorem}{Theorem}
\newtheorem{lemma}[theorem]{Lemma}
\newtheorem{claim}[theorem]{Claim}
\newtheorem{corollary}[theorem]{Corollary}
\newtheorem{prop}[theorem]{Proposition}
\newtheorem{note}{Note}
\newtheorem{definition}[theorem]{Definition}
\newtheorem{example}[theorem]{Example}
\numberwithin{equation}{section}



\title{The symmetric groups $S_n, n\geq 4$, and finite non-abelian 
simple groups are not embeddable in any Riordan group}


\author[1]{Tian-Xiao He}
\ead{the@iwu.edu}
\author[2]{Nikolai A. Krylov\fnref{fn2}\corref{cor2}}
\ead{nkrylov@siena.edu}
\affiliation[1]{organization={Department of Mathematics, Illinois Wesleyan University},   
		addressline={1312 Park Street}, 
            	city={Bloomington},
            	postcode={61702}, 
            	state={IL},
            	country={USA}}
\affiliation[2]{organization={Department of Mathematics, Siena University},   
		addressline={515 Loudon Road}, 
            	city={Loudonville},
            	postcode={12211}, 
            	state={NY},
            	country={USA}}

\cortext[cor2]{Corresponding author}
\fntext[fn2]{Nikolai Krylov is supported by an AMS-Simons 
Research Enhancement Grant for Primarily Undergraduate Institution Faculty.}

\begin{abstract}
We prove that the symmetric group of degree greater than three cannot be embedded
into the Riordan group with coefficients in any commutative ring. We also prove the 
impossibility to embed finite non-abelian simple groups. As a closely related topic, 
we show why all truncated Riordan groups are solvable, in stark contrast to the 
unsolvability of the infinite-sized Riordan groups. Finally, we give an explicit 
embedding of the alternating group $A_4$ into the Lagrange subgroup with 
coefficients in a certain commutative ring, and prove that $A_4$ cannot be 
embedded into a substitution group.
\end{abstract}


\end{frontmatter}




\section{Introduction}
\label{sec1}



This article is a continuation of our paper \cite{HeKrylov}, where we started a discussion about 
what groups can be represented by Riordan arrays. We proved in particular, that the
symmetric group of degree three $S_3$ has no faithful representation as a subgroup of the 
Riordan group over the complex numbers, but can be embedded as a subgroup of the 
Riordan group over a field of characteristic three. We also asked questions about 
existence of such representations for $S_n, n\geq 4$ over other fields of finite 
characteristic, or more generally, over an arbitrary commutative ring with identity 
(see Question 3. and Question 4. at the end of \cite{HeKrylov}). 
Here we answer those questions in the negative, the symmetric group $S_n, n\geq 4$ 
is not embeddable in any Riordan group.

The Riordan group, introduced by Shapiro et al. \cite{Shapiro1}, contains infinite lower 
triangular matrices called Riordan arrays, whose columns consist of the coefficients of 
certain formal power series. For the detailed introduction to the subject, the reader is 
urged to study the books by Barry \cite{Barry} and Shapiro et al. \cite{Shapiro2}, 
and a survey article by Davenport et al. \cite{Davenport}. Here are a few basic 
definitions and notations, which we will be using in the text. 

Let us fix an arbitrary commutative ring with identity, say $\mathbb D$, which may 
have zero divisors, and use $\mathbb D^*$ to denote the multiplicative group of its 
units. The set of all formal power series (f.p.s.) in indeterminate $t$ with coefficients 
in $\mathbb D$ is denoted by $\F = {\mathbb D}[\![$$t$$]\!]$. The \emph{order} of 
$f(t)  \in \F$, $f(t) =\sum_{k=0}^\infty f_kt^k$ ($f_k\in {\mathbb D}$), is the minimal 
number $r\in{\mathbb N_0=\{0\}\cup\natu}$ such that $f_r \neq 0$, and the set of 
formal power series of order $r$ is denoted by ${\F}_r$. Let $g(t) \in \F_0$ and 
$f(t) \in \F_1$; the pair $\bigl(g(t) ,\,f(t)\bigr)$ defines the {\em Riordan array} 
$$
A =(d_{n,k})_{n,k\geq 0}=\bigl(g(t) ,\,f(t)\bigr)
$$ 
having
\begin{equation}\label{1}
d_{n,k} = [t^n]g(t) f(t)^k,
\end{equation} 
where $[t^n]h(t)$ denotes the coefficient of $t^n$ in the expansion of a f.p.s. $h(t)$. 
When $\mathbb D$ is a field, the set of all such Riordan arrays forms a multiplicative group, 
called the {\sl Riordan group}, and denoted by ${\cal R}({\mathbb D})$. The group operation 
$*$ is written in terms of the f.p.s. as 
$$
\bigl(g_1(t) ,\,f_1(t)\bigr) * \bigl(g_2(t) ,\,f_2(t)\bigr) = \bigl(g_1(t) g_2(f_1(t)),\,f_2(f_1(t))\bigr),
$$
(which will be omitted from now on) with the Riordan array $I = (1,\,t)$ acting 
as the group identity. The inverse of the Riordan array $\bigl(g(t) ,\,f(t)\bigr)$ is the pair
$$
\bigl(g(t) ,\,f(t)\bigr)^{-1} = \left(\frac{1}{g(\bar{f}(t))},\,\bar{f}(t)\right),
$$
where we used the standard notation $\bar{f}(t)$ for the compositional inverse of $f(t)$. Thus, 
$\bar{f}(f(t)) = t$ and $f(\bar{f}(t)) = t$. If $\mathbb D$ is not a field, but only a commutative ring 
with identity, to guarantee the existence of the inverse $\bigl(g(t) ,\,f(t)\bigr)^{-1}$, the coefficients 
$g_0$, and $f_1$ must be units, i.e. we will require in addition, that $g_0,f_1\in \mathbb D^*$. 

The Riordan group ${\R}({\mathbb D})$, is the semidirect product of two 
proper subgroups: the Appell subgroup and the Lagrange (or associated) 
subgroup (\cite{Barry}, \cite{Shapiro2}). The Appell subgroup 
${\cal A}({\mathbb D})$ is abelian, normal, and consists of the 
Riordan arrays $\bigl(g(t) ,\,t\bigr)$. 

The Lagrange subgroup ${\cal L}({\mathbb D})$ consists of the Riordan arrays 
$\bigl(1 ,\,f(t)\bigr)$, and contains the {\sl substitution group of formal power 
series} $\cal{J}({\mathbb D})$, which was introduced in the 1950s (see papers 
by Jennings \cite{Jennings}, Johnson \cite{Johnson}, and Babenko \cite{Babenko}). 
Elements of the substitution group, as Riordan arrays, are the 
pairs $\bigl(1 ,\,f(t)\bigr)$, where the f.p.s. $f(t) =\sum_{k \geq 1} f_kt^k$ 
has the first coefficient $f_1 = 1$. Hsu, Shiue and one of the 
authors \cite{HHS} analyzes the 
Sheffer group - a version of the substitution (Sheffer) group - as it relates to, 
and sits within, Riordan subgroups, elucidating how substitution-type 
operations embed within the broader Riordan framework. Furthermore, it is 
proved in Proposition 5. of \cite{Johnson} that if the multiplicative identity 1 has 
infinite order in the additive group ${\mathbb D}^+$, the substitution group 
${\cal J}({\mathbb D})$ contains a copy of $F_2$, a free group of rank two. 
As a subgroup of the Riordan group, $F_2$ is generated by the elements 
$$
\bigl(1,\, t/(1 + 3t)\bigr) ~~~ \mbox{and} ~~~ 
\left(1,\, t/\bigl(\sqrt[3]{{1 + (3t)^3}}\bigr)\right).
$$
The alternating group $A_5$ is generated by two 
permutations, $(12345)$ and $(12)(34)$, so naturally is a factor group of 
$F_2$. Since $A_5$ is not solvable, $F_2$ is not solvable, and hence 
neither of the groups 
$$
{\cal J}({\mathbb D}) \leq {\cal L}({\mathbb D}) \leq {\R}({\mathbb D})
$$ 
is solvable. If we take the coefficient ring ${\mathbb D}$ 
to be a finite field on $p$ elements ${\mathbb D} = \inte_p = \inte/p\inte$, 
the substitution group ${\cal J}(\inte_p)$ is also known under the name of the 
{\sl Nottingham group} (see Camina's survey article \cite{Camina2} for 
more details). This Nottingham group is a pro-$p$ group, i.e. it is the 
inverse limit of a tower of finite $p$-groups, and if $p > 2$, it is known to 
be a finitely presented with two generators (see Ershov \cite{Ershov}). 
Camina proved in \cite{Camina1} that every finitely generated pro-$p$ group 
can be embedded, as a closed subgroup, in ${\cal J}(\inte_p)$.
Since the pro-$p$ completion $\widehat{F}_2$ of the free group 
$F_2 \leq \widehat{F}_2$ is a countably based pro-$p$ group, it implies that $F_2$ 
can be embedded in ${\cal J}(\inte_p)$, and hence ${\cal J}(\inte_p)$ 
and ${\R}(\inte_p)$ are not solvable. Moreover, Szegedy proved that 
two randomly chosen elements of the Nottingham group generate a free 
subgroup with probability 1, and that ${\cal J}(\inte_p)$ contains a dense 
free subgroup of rank two (see Theorems 1. and 2. in \cite{Szegedy}).
On the other hand, the symmetric group of degree four is solvable, and 
\begin{equation}
\label{Klein}
\{e\} \lhd K_4 = \{e,(12)(34),(13)(24),(14)(23)\} \lhd A_4 \lhd S_4.
\end{equation}
is a subnormal solvable series, where $A_4$ is the alternating group, 
and $K_4$ is the Klein four-group. Our main goal here is to show 
that there exists no embedding of $S_4$ (and therefore of $S_n, n > 4$) 
into the Riordan group ${\R}({\mathbb D})$ 
with coefficients in any commutative ring ${\mathbb D}$. 

From now on, we assume that $\mathbb D$ is a fixed, arbitrary commutative 
ring with identity 1 (unless otherwise stated), 
and write the Riordan group ${\R}({\mathbb D})$, the Appell subgroup 
${\cal A}({\mathbb D})$, the Lagrange subgroup ${\cal L}({\mathbb D})$, 
and the substitution group ${\cal J}({\mathbb D})$ with coefficients in 
$\mathbb D$ respectively as ${\R}, {\cal A}, {\cal L}$, and ${\cal J}$.

Our interest in studying embeddings of the symmetric group $S_n$ into the 
Riordan group ${\R}$ is motivated by the fundamental importance of these 
groups particularly in combinatorics. The symmetric group $S_n$ is the 
group of all bijective rearrangements of $n$ objects with the operation of 
composition of such rearrangements. It is called symmetric because 
symmetry of a set can be viewed as a structure preserving map of the 
set into itself. Symmetric groups model the idea of symmetry in the 
most general form, and any symmetry (in geometry, algebra, or physics) is often
related to a permutation of objects. Symmetric groups arise naturally 
almost everywhere. The group of symmetries of a regular $n$-gon 
is a subgroup of $S_{2n}$. The Galois group of a polynomial $p(x)$ is a 
subgroup of all permutations of the roots of $p(x)$. Linear representations 
of symmetric groups utilize matrix algebra in studying 
symmetries of objects with crucial applications in physics modeling 
particle exchange symmetry (e.g. bosons vs fermions), and in chemistry. 
Furthermore, in computer science, the symmetric group provides a framework for 
finding optimal orderings, sorting algorithms, and analyzing network optimization.
In survey, the symmetric group $S_n$ is a central object in mathematics because 
it captures the essence of permutation, symmetry, and structure. It is both a tool and 
a model that expands and interweaves many mathematical themes, and appears 
almost in every branch of mathematics.

The Riordan arrays and Riordan group, introduced less than four decades ago, 
already proved to play not only a unifying role in enumerative combinatorics, 
but to inspire new research directions in such fields of mathematics, as 
algebra, algebraic and enumerative combinatorics, generating functions, 
special functions, posets, lattice paths, orthogonal and Sheffer 
polynomials, probability, group theory, graph theory, topology, etc. 

We see a significant overlap in the application areas of symmetric groups and Riordan 
groups, and it is natural and tempting to study if there is a deeper interconnection  
between them. 

In many problems, $S_n$ and the Riordan group represent different aspects 
of the problem. For instance, the Riordan array $\bigl(C(t), \, tC(t)\bigr)$, where 
$C(t)=(1-\sqrt{1-4t})/2t$ is the Catalan generating function, can be used to model 
the structure of $132$-avoiding permutations (cf. B\'ona \cite{Bon}), because 
of the following reason: The entries $d_{n,k}$ of $\bigl(C(t),\, tC(t)\bigr)$ 
count the number of $132$-avoiding permutations of size $n$ that can be 
decomposed into $k+1$ blocks, where each block is itself a $132$-avoiding 
permutation, and the entire permutation respects a hierarchical (non-crossing) 
structure. Another important motivation for this study may be the property of 
symmetric groups as a ``universal container" for all finite group structures, which 
is proposed by Cayley's theorem: every finite group can be embedded in some $S_n$.
In this paper we will show that Cayley's theorem no longer holds for groups 
such as ${\cal R}$ and ${\cal R}_n$, where the latter is an $n$-truncation of 
${\cal R}$ (a Riordan group consisting of $(n+1)\times (n+1)$ size Riordan arrays).

The structure of this article is the following. Truncated Riordan 
groups ${\cal R}_n$ and truncated Appell and Lagrange subgroups,  
${\cal A}_n$ and ${\cal L}_n$, 
will be discussed in the next section. In particular, we will explain why all 
groups ${\cal R}_n$ and ${\cal L}_n$ are solvable. In section $3$, 
we will use these truncated groups to prove that neither the symmetric 
group $S_4$ nor any finite non-abelian simple group 
can be embedded as a subgroup of the Riordan 
group over any commutative ring. At the end we will show how 
one can embed the alternating group $A_4$ into the 
Lagrange subgroup, and prove that there exists no embedding 
of $A_4$ into the substitution group of f.p.s., and in particular, into 
the Nottingham group ${\cal J}({\mathbb F}_p)$.


\section{Truncated Riordan groups}

Each Riordan array naturally reduces to a lower-triangular square matrix 
of size $k \times k$ by taking the first $k$ rows and columns. This idea 
leads to the description of the Riordan group $\R$ as the inverse limit 
of an inverse sequence of groups of finite size matrices, which was presented 
in \cite{Luzon16}. It was further exploited in \cite{Luzon17} to give a formula 
for all Riordan involutions, and also in \cite{Luzon17B}, 
where the truncated groups ${\cal R}_n(\real)$ and ${\cal R}_n(\comp)$ 
were given the structure of a differentiable manifold (complex respectively).
Let us briefly mention the construction, because a part of it plays a 
key role in our proofs below. For further details we refer the 
reader to \S3 of \cite{Luzon16}, and \S2 of \cite{Luzon17}, where all the Riordan 
groups are considered over a field of characteristic zero. The domain extension to 
an arbitrary ring is straightforward and we leave the formal details to the reader. 

Take an arbitrary $n\in\natu$, and consider the general linear group of all invertible 
$(n+1)\times (n+1)$ matrices ${\rm GL}(n+1,{\mathbb D})$ with coefficients in the 
fixed ring $\mathbb D$. There is a natural {\sl truncation} homomorphism
\begin{equation}
\label{Pin}
\Pi_n: {\R} \longrightarrow {\rm GL}( n+1, {\mathbb D})
\end{equation}
defined by
$$
\Pi_n\bigl((d_{i,j})_{i,j\geq 0}\bigr) = (d_{i,j})_{0\leq i,j\leq n}.
$$
The image of this homomorphism is a subgroup of ${\rm GL}( n+1, {\mathbb D})$, 
which is denoted by ${\R}_n = \Pi_n({\R})$. Clearly, we can obtain ${\R}_n$ from 
${\R}_{n+1}$ by deleting the last row and column. This operation gives 
another natural homomorphism 
\begin{equation}
\label{FinPin}
P_n: {\R}_{n+1} \longrightarrow {\R}_{n}
\end{equation}
formally defined by
$$
P_n\bigl((d_{i,j})_{0\leq i,j\leq n+1}\bigr) = (d_{i,j})_{0\leq i,j\leq n}
$$
(see Definition 3. in \cite{Luzon16}). Homomorphisms $\Pi_{n+1},\Pi_n$ and $P_n$ are 
connected in the commutative diagram
\begin{equation}
\label{ComD}
\begin{tikzcd}
 & {\R} \arrow[dl,"\Pi_{n+1~}" left, near start] \arrow[dr, "\Pi_n"] & \\
{\R}_{n+1} \arrow{rr}{P_n} & & {\R}_n\\
\end{tikzcd}
\end{equation}
that is $\Pi_n = P_n\circ\Pi_{n+1}$, and the Riordan group ${\R}$ is isomorphic to the 
inverse limit $\lim\limits_{\longleftarrow}\{({\R}_n)_{n\in\natu},(P_n)_{n\in\natu}\}$ 
(see Proposition 4. in \cite{Luzon16}).

The homomorphism $P_n$ is onto for all $n\in\natu$, i.e. ${\rm Im}(P_n) = {\R}_n$. 
The kernel of $P_n$ is given implicitly in \cite{Luzon17} (see Proposition 2. and the 
discussion right above it). Since this kernel plays the fundamental role in our proofs, 
we give an explicit statement with a detailed proof in the following

\begin{lemma}
The kernel of $P_0$ is isomorphic to the semidirect product 
$$
\ker(P_0) \cong {\mathbb D} \rtimes_{\varphi} {\mathbb D}^*,
$$
where $\varphi: {\mathbb D}^*\to Aut({\mathbb D})$ is a homomorphism 
defined by $\varphi(a)(b) = \varphi_a(b) = ab$.
The kernel of $P_n, n\geq 1$ is isomorphic to the direct product 
$$
\ker(P_n) \cong {\mathbb D} \times {\mathbb D}.
$$
\end{lemma}
\begin{proof}
The statement about the kernel of $P_0$ is proved in Lemma 1 of \cite{HeKrylov}, so 
assume that $n\geq 1$. The main observation here is that every entry  of a Riordan array
$d_{i,j}$, where $i,j\geq 1$, is determined be the A-sequence and the entries in the
row above. Therefore, if $P_n(M) = I \in{\R}_n$ for some Riordan array $M$, 
we must have $a_0 = 1$ and $a_i = 0$ for all $i\in\{1,\ldots,n - 1\}$, where 
$A(t) = \sum\limits_{i\geq 0} a_it^i$ is the generating function of the 
A-sequence of $M$. Hence, for all $j\in \{2,\ldots,n\}$ we also have
$$
d_{n+1,j} = \sum\limits_{k = 0}^{n-j} a_k\cdot d_{n,j-1+k}  = 0,
$$
and of course, $d_{n+1,n+1} = d_{n,n} = 1$. On the other hand, $d_{n+1,0}$ can be an arbitrary 
element of the additive abelian group $\mathbb D$, and $d_{n+1,1}$ can also be any 
element of $\mathbb D$ because $d_{n+1,1} = a_n\cdot d_{n,n} = a_n$, and the truncated 
Riordan arrays from ${\R}_n$ impose no restrictions on the corresponding terms $a_n$ of their 
A-sequences. Next, take arbitrary $\alpha_1,\alpha_2,\beta_1,\beta_2\in \mathbb D$, 
and consider the f.p.s. 
$$
g_i(t) = 1 + \alpha_i t^{n+1} + \sum\limits_{k \geq n+2} 0t^k ~~~ 
\mbox{and} ~~~ f_i(t) = t + \beta_it^{n+1}  + \sum\limits_{k \geq n+2} 0t^k,
$$
where $i\in \{1,2\}$. It is clear that for each $i\in\{1,2\}$, the Riordan array $\bigl(g_i,f_i\bigr)$  
belongs to the kernel of $\Pi_n$. Moreover, the product Riordan array equals
$$
(g_1,f_1) (g_2,f_2) = 
\bigl(1 + \alpha_1t^{n+1},t+\beta_1t^{n+1}\bigr)\bigl(1 + \alpha_2t^{n+1},t+\beta_2t^{n+1}\bigr)
$$ 
$$
= 
\bigl(1 + (\alpha_1+\alpha_2)t^{n+1} + \ldots h.o.t., ~ 
t + (\beta_1+\beta_2)t^{n+1} + \ldots h.o.t.\bigr),
$$
where $h.o.t.$ stands for higher-order terms.
Let us denote this product matrix by $\bigl(m_{i,j}\bigr)_{i,j\geq 0}$. 
Then, the nonzero entries of the truncated matrix 
$$
\Pi_{n+1}\left(\bigl(m_{i,j}\bigr)_{i,j\geq 0}\right) = 
\bigl(m_{i,j}\bigr)_{0\leq i,j\leq n+1}
$$ 
will be
$$
m_{k,k} = 1 ~ \mbox{for} ~ 0\leq k\leq n+1, ~~~ \mbox{and} ~~~ 
m_{n+1,0} = \alpha_1 + \alpha_2, ~~ m_{n+1,1} = \beta_1 + \beta_2.
$$
It implies that the matrices $\Pi_{n+1}\bigl((g_1,f_1)\bigr)$ and 
$ \Pi_{n+1}\bigl((g_2,f_2)\bigr)$ commute 
$$
\Pi_{n+1}\bigl((g_1,f_1)\bigr) \cdot \Pi_{n+1}\bigl((g_2,f_2)\bigr) = 
\Pi_{n+1}\bigl((g_2,f_2)\bigr) \cdot \Pi_{n+1}\bigl((g_1,f_1)\bigr),
$$
and therefore $\ker(P_n) \cong {\mathbb D} \times {\mathbb D}$, as required.
\end{proof}
We can rewrite this lemma in the form of a short exact sequence  
\begin{equation}
\label{ShExSeq}
1\longrightarrow {\mathbb D} \times {\mathbb D} \stackrel{\l_{n+1}}\longrightarrow {\R}_{n+1} 
\stackrel{P_n}\longrightarrow
{\R}_n\longrightarrow 1, ~~\forall n\geq 1
\end{equation}
where the homomorphism $l_{n+1}$ is defined for a pair 
$(\alpha,\beta)\in {\mathbb D} \times {\mathbb D}$ by 
$$
l_{n+1}\bigl((\alpha,\beta)\bigr) = \Pi_{n+1}\bigl((1+\alpha t^{n+1}, t + \beta t^{n+1})\bigr) \in {\R}_{n+1}.
$$

Next we show that if a finite group $G$ is embeddable into the Riordan group $\R$, 
then there is a truncated Riordan group ${\R}_n$ (over the same ring $\mathbb D$, of 
course), which contains an isomorphic copy of $G$ as a subgroup.

\begin{prop}
Let $G$ be a finite group and $\mu:G\to {\R}$ a monomorphism. Then there 
exists $n\in\natu$ such that the composition 
$$
\Pi_n\circ\mu: G \longrightarrow {\R} \longrightarrow {\R}_n
$$
is a monomorphism as well.
\end{prop}
\begin{proof}
Assume that $G$ is not a trivial group, and take any nontrivial element $g\in G$. Since $\mu(g)$ 
is not the identity Riordan array, $\exists m_g\in\natu$ such that 
$$
\Pi_{m_g}\bigl(\mu(g)\bigr) \neq I\in {\R}_{m_g}.
$$
Take the largest of all such numbers
$$
n = \max\{m_g~|~g\in G\} \in\natu,
$$
which does exist since $|G|<\infty$, and consider the composition 
$$
\Pi_n\circ\mu: G \longrightarrow {\R} \longrightarrow {\R}_n.
$$
Since $\Pi_{n}\bigl(\mu(g)\bigr) \neq I\in {\R}_n$ for every nontrivial $g\in G$, the composition 
$\Pi_n\circ \mu$ is a monomorphism. 
\end{proof}

\begin{note} 
It is easy to see that the statement above is false when $G$ has infinitely 
many elements. As a simple counterexample, consider the Riordan group with 
coefficients in an arbitrary finite commutative ring with identity, and the free 
abelian group $G = \langle (1 + t, \, t)\rangle\hookrightarrow {\R}$
generated by the element $(1 + t, \, t)$. Since the ground ring is finite, any truncated 
Riordan group ${\R}_n$ will have only finitely many elements, and can not contain 
$G$ as a subgroup. For a counterexample with an infinite ground ring, we can take 
the Riordan group over the polynomial ring ${\mathbb D}  = \inte_2[x]$, and consider 
the same subgroup $G = \langle (1 + t,\, t)\rangle$. Since in this ring 
$\bigl((1 + t)^{2^k}, \,t\bigr) = \bigl(1 + t^{2^k},\, t\bigr) \neq e$, it is clear that for any $n\in\natu$, 
$\Pi_n\bigl(\mu\bigl((1 + t)^{2^k},\, t\bigr)\bigr) = I \in {\cal R}_n$ when $k > \ln(n)/\ln(2)$.
\end{note}
\begin{note}
We also would like to notice here that restricting the truncation homomorphism 
$\Pi_n,\forall n\in \natu$ onto any subgroup $G\leq {\R}$ will produce a homomorphism
$$
\left.\Pi_n\right|_{G} : G \lra  {\rm GL}( n+1, {\mathbb D}).
$$
In particular, we can restrict $\Pi_n$ onto the Appell and Lagrange subgroups 
of ${\R}$, to obtain the corresponding truncated subgroups 
$$
{\cal A}_n := \left.\Pi_n\right|_{\cal A} ~~~ \mbox{and} ~~~  
{\cal L}_n := \left.\Pi_n\right|_{\cal L}.
$$
Since the Riordan group ${\R}$ is the semidirect product 
${\R}\cong {\cal A} \ltimes {\cal L}$, it is a straightforward exercise to show 
that for each $n\in\natu$ we have the semidirect product of the 
truncated groups ${\R}_n\cong {\cal A}_n \ltimes {\cal L}_n$ as well (c.f. Proposition 32. in 
\cite{Luzon17B} when the ground ring is a field). Deleting the last row and column 
in the matrix representing an element of ${\cal A}_{n+1}$ or ${\cal L}_{n+1}$ produces 
correspondingly an element of ${\cal A}_{n}$ or ${\cal L}_{n}$. Thus we obtain 
homomorphisms, which are the restrictions of the homomorphism $P_n$ onto the 
subgroups ${\cal A}_{n+1}$ and ${\cal L}_{n+1}$ respectively. 
\end{note}

Our next Lemma can be considered as a direct corollary of Lemma 1.

\begin{lemma} For the truncated Appell and Lagrange subgroups over the ring 
$\mathbb D$ we have 
$$
{\cal A}_0 \cong {\mathbb D}^* ~~~ \mbox{and} ~~~ 
{\cal L}_0 \cong \{e\},~{\cal L}_1 \cong {\mathbb D}^*.
$$
The kernel $\ker\left(\left.P_0\right|_{{\cal A}_1}\right) \cong {\mathbb D}$, and if $n\geq 1$, 
the kernels of $\left.P_n\right|_{{\cal A}_{n+1}}$ and $\left.P_n\right|_{{\cal L}_{n+1}}$ 
are isomorphic to the additive group of the ring $\mathbb D$ as well. 
In other words, for all positive $n\in\natu$ we have 
two short exact sequences (compare with \ref{ShExSeq})
\begin{equation}
\label{ShExSeq2A}
1\longrightarrow {\mathbb D} \stackrel{\nu_{n+1}}\longrightarrow {\cal A}_{n+1} 
\stackrel{\left.P_n\right|_{{\cal A}_{n+1}}}\longrightarrow
{\cal A}_n\longrightarrow 1,
\end{equation}
\begin{equation}
\label{ShExSeq2L}
1\longrightarrow {\mathbb D} \stackrel{\mu_{n+1}}\longrightarrow {\L}_{n+1} 
\stackrel{\left.P_n\right|_{{\cal L}_{n+1}}}\longrightarrow
{\L}_n\longrightarrow 1,
\end{equation}
where for any $d\in\mathbb D$, $\nu_{n+1}(d): = \bigl(1+dt^{n+1},\,t\bigr)$, 
and $\mu_{n+1}(d): = \bigl(1,\,t + dt^{n}\bigr)$.
\end{lemma}
\begin{proof}
(\ref{ShExSeq2A}) follows immediately from the formula for the product of two 
polynomials $[t^n]\bigl((1+d_1t^n)(1+d_2t^n)\bigr)=d_1+d_2$, and 
(\ref{ShExSeq2L}) follows from the formula for the composition of two polynomials
$f(t) = t + f_nt^n$ and $h(t) = t + h_nt^n$, since 
$f(h(t)) = t + h_nt^n + f_n(t + h_nt^n)^n$, and $[t^n]f(h(t)) = h_n + f_n$. 
Further details are left to the reader.
\end{proof}

Note that neither of  (\ref{ShExSeq2A}) or (\ref{ShExSeq2L}) splits in general. For 
example, if $\mathbb D = {\mathbb F}_2$ ${\cal A}_2\cong C_4$, the cyclic group of 
order 4, which is generated by the array $(1+t,\,t)$.
The results of Lemma 1. and Lemma 3. can be visualized and put together into one 
commutative diagram, where all horizontal and vertical sequences are short exact 
sequences. We do not name the corresponding homomorphisms in the diagram to 
keep the notations simple, but we explained all the homomorphisms above 
(c.f. (\ref{ShExSeq}), (\ref{ShExSeq2A}) and (\ref{ShExSeq2L})). We present here the 
case when $n\geq 1$, the case of $n = 0$ is similar. Also notice that all horizontal 
sequences split, since the groups in the middle are the semidirect products. 
\begin{equation}
\label{BigComD}
\begin{CD}
@.  1 @. 1  @. 1 @.\\
@.  @ VVV  @ VVV  @ VVV @.\\
1 @>>> {\mathbb D} @>>> {\mathbb D}\times {\mathbb D} @>>> {\mathbb D} @>>> 1\\
@.  @ VVV  @ VVV  @ VVV @.\\
1 @>>> {\cal A}_{n+1} @>>> {\cal R}_{n+1} @>>> {\cal L}_{n+1} @>>> 1\\
@.  @ VVV  @ VVV  @ VVV @.\\
1 @>>> {\cal A}_n @>>> {\cal R}_n @>>> {\cal L}_n @>>> 1\\
@.  @ VVV  @ VVV  @ VVV @.\\
@.  1 @. 1  @. 1 @.\\
\end{CD}
\end{equation}

\begin{note}
Since a group $G$ is solvable if and only if for a normal subgroup $K$ of $G$, both 
$K$ and $G/K$ are solvable, and ${\mathbb D}^* \cong {\cal R}_0\cong {\cal L}_1$ 
is abelian, the diagram (\ref{BigComD}) together with the induction imply that the truncated 
groups ${\cal R}_n$ and ${\cal L}_n$ are solvable for all $n\in \natu$. See Remarks 
14 and 15 at the end of \cite{Luzon23}, and also the first paragraph in section 
three of \cite{Luzon17B} for different arguments explaining the solvability of 
${\cal R}_n$ over a field. Recall also the introduction section, where we explained 
why the Riordan group ${\cal R}$ is not solvable in general.
\end{note}


\section{Impossibility of embedding $S_n\to {\R},n\geq 4$}

In this section we prove that the symmetric group $S_n$ of degree $n\geq 4$ 
can not be embedded as a subgroup of the Riordan group with coefficients in 
any commutative ring. Our proof is based on the fact that $S_4$ has exactly 
four normal subgroups $\{\{e\}, K_4, A_4, S_4\}$, (recall (\ref{Klein})). 
We will use Lemma 1. to show that $S_4$ can not be embedded into ${\R}$, 
and hence, none of the groups $S_n,~n\geq 5$ can. As above, we 
assume that the Riordan group $\R$ has coefficients in a fixed commutative 
ring $\mathbb D$ with identity, denoted by 1. 

\begin{theorem}
If $n\geq 4$, there exists no monomorphism 
\begin{equation}
\label{MonomSn}
\mu: S_n\to {\R}.
\end{equation}
\end{theorem}
\begin{proof}
Assume to the contrary that there is a monomorphism 
$\mu: S_4\to {\R}$. Then according to Proposition 2, $\exists N\in\natu$, 
such that the composition of $\mu$ with $\Pi_N$ will be a monomorphism 
$\Pi_N\circ \mu : S_4\to {\R}_N$. Let $m$ be the smallest of all such 
possible numbers $N$. To simplify the notations we will denote the 
composition $\Pi_m\circ \mu$ simply by $\mu:S_4\to {\R}_m$ in this proof.
Since ${\R}_0\cong {\mathbb D}^*$ is abelian and can not contain $S_4$, 
we have a monomorphism $\mu: S_4 \to {\R}_m,$ where $m\geq 1$. 
Suppose first that $m = 1$, then we have
$$
\ker(P_0\circ\mu) \lhd S_4 ~~~ \mbox{and} ~~~ \ker(P_0\circ\mu)\leq  
{\mathbb D} \rtimes_{\varphi} {\mathbb D}^*.
$$
Since ${\R}_0$ is commutative, it can not contain a copy of $S_4$ or 
$S_3\cong S_4/K_4$. Hence the kernel $\ker(P_0\circ\mu)$ can not  
be trivial or $K_4$ respectively.  By a similar reason we can not have 
$\ker(P_0\circ\mu) = S_4 \leq {\mathbb D} \rtimes_{\varphi} {\mathbb D}^*$. 
Indeed, since $\{\{e\}, K_4, A_4, S_4\}$ are the only normal subgroups of $S_4$, 
in such a case the short exact sequence 
\begin{equation}
\label{ExSeq2}
1\longrightarrow {\mathbb D} \longrightarrow {\mathbb D} \rtimes_{\varphi} {\mathbb D}^*
\longrightarrow {\mathbb D}^* \longrightarrow 1,
\end{equation}
would imply that either ${\mathbb D}$ or ${\mathbb D}^*$ 
contains a nonabelian subgroup. Thus, we end up with the last possible option 
when $\ker(P_0\circ\mu) = A_4$. 
This option can be described by the following commutative diagram, where all 
arrows represent group homomorphisms, and horizontally we have two short 
exact sequences.
\begin{equation}
\label{ComDiag2}
\begin{CD}
1 @>>> A_4  @>>> S_4 @>>> \inte_2 @>>> 1\\
@. @VV{\left.\mu \right|_{A_4}}V @ VV{\mu}V @ VVV @.\\
1 @>>> {\mathbb D} \rtimes_{\varphi} {\mathbb D}^* @>>>  {\R}_1 @>P_0>> {\R}_0  @>>> 1
\end{CD}
\end{equation}
Homomorphism $\mu$ is a monomorphism by the assumption, the restriction of $\mu$ 
on $A_4$, which we denote by $\left.\mu \right|_{A_4}$, is also a monomorphism. 
A standard diagram chasing argument (left to the reader) shows that the last 
vertical arrow in (\ref{ComDiag2}) is a monomorphism as well. If it happens that 
${\R}_0 \cong {\mathbb D}^*$ does not contain a cyclic group of order 2, the proof 
stops here. Otherwise, the diagram (\ref{ComDiag2}) implies that all elements of 
$S_4$ can be represented by $2\times 2$ lower triangular matrices of the form
$$
\begin{pmatrix}
a & 0\\
b & ac\\
\end{pmatrix}, 
$$
where $a\in \inte_2\leq {\mathbb D}^*, b\in {\mathbb D}, c\in{\mathbb D}^*$. 
Moreover, the second short exact sequence in (\ref{ComDiag2}) actually splits, 
i.e. there is a homomorphism $s:{\R}_0 \longrightarrow {\R}_1$ 
such that $P_0\circ s = Id_{{\R}_0}$. Indeed, define this homomorphism $s$ by 
$s(a) : = \Pi_1\Bigl((a,t)\Bigr)$, or 
$$
s(a) : = \begin{pmatrix}
a & 0\\
0 & a\\
\end{pmatrix}, ~ \forall a\in {\mathbb D}^*.
$$
Since for all $a\in {\mathbb D}^*$ and $M\in {\R}_1$, $s(a)$ and $M$ commute, 
the group ${\R}_1$ is isomorphic to the direct product 
\begin{equation}
\label{Isom1}
{\R}_1\cong \left({\mathbb D} \rtimes_{\varphi} {\mathbb D}^*\right) \times {\mathbb D}^*.
\end{equation}
An isomorphism (\ref{Isom1}) together with the diagram (\ref{ComDiag2}) imply 
that $S_4$ has a transposition, which commutes with every permutation from the 
alternating subgroup $A_4$. This is a contradiction, since $S_4$ is not isomorphic 
to the direct product $A_4\times \inte_2$, and therefore we can not have 
$\ker(P_0\circ\mu) = A_4$.

Now let us suppose that there is a monomorphism $\mu: S_4\longrightarrow {\R}_m$, 
where the smallest such $m\geq 2$. As follows from the second part of Lemma 1, 
in such a case we must have 
$$
\ker(P_{m-1} \circ\mu) \lhd S_4 ~~~ \mbox{and} ~~~ \ker(P_{m-1} \circ\mu)\leq  
{\mathbb D} \times {\mathbb D}.
$$
Minimality of $m\geq 2$, and commutativity of ${\mathbb D} \times {\mathbb D}$ 
rule out for $\ker(P_{m-1} \circ\mu)$ the options \{\{$e$\}, $A_4$, $S_4$\}, and  
we need to discuss only one possibility, if $\ker(P_{m-1} \circ\mu) = K_4$. 
This option has the following description in terms of the commutative diagram
\begin{equation}
\label{ComDiag3}
\begin{CD}
1 @>>> K_4  @>>> S_4 @>>> S_3 @>>> 1\\
@. @VV{\left.\mu \right |_{K_4}}V @ VV{\mu}V @ VVV @.\\
1 @>>> {\mathbb D} \times {\mathbb D} @>l_{m}>>  {\R}_m @>P_{m-1} >> {\R}_{m-1}  @>>> 1
\end{CD}
\end{equation}
where, as in (\ref{ComDiag2}), horizontally we have two short exact sequences 
(recall (\ref{ShExSeq}) for the definition of $l_m$), and all vertical arrows represent 
monomorphisms. The idea here is to consider the product of a 3-cycle from $\mu(S_4)$ 
(there are 8 elements of order 3 in $S_4$) with an involution from $\mu(K_4)$. 
We show that the diagram (\ref{ComDiag3}) implies that such product can not have  
order 3, which contradicts the structure of $S_4$, where the 
product of any 3-cycle with any involution from $K_4$ given by (\ref{Klein}), is 
again a 3-cycle. We will need the following 

\begin{claim}
Let $(g,f)\in \mu(S_4) \subseteq {\R}$ be a Riordan array of order 3 
with the f.p.s. expansions $g(t) = \sum\limits_{i\geq 0} g_it^i$ and 
$f(t) = \sum\limits_{i\geq 1} f_it^i$. Then $g_0 = f_1 = 1$.
\end{claim}
\begin{proof}[Proof of the claim]
Every 3-cycle $(abc)=(ab)(bc)$ is a product of two involutions. 
Hence, there are involutions in $\R$, say 
$$
(u,v) = \left(\sum\limits_{i\geq 0} u_it^i,\sum\limits_{i\geq 1} v_it^i\right),~~ 
\mbox{and} ~~ (p,q) = \left(\sum\limits_{i\geq 0} p_it^i,\sum\limits_{i\geq 1} q_it^i\right)
$$
such that 
$$
\bigl(g(t),f(t)\bigl) = \bigl(u(t),v(t)\bigl)\bigl(p(t),q(t)\bigl) = \bigl(u(t)p(v(t)),q(v(t))\bigr).
$$
It means, in particular, that
\begin{equation}
\label{formula1}
u_0p_0 = g_0, ~~~ \mbox{and} ~~~ q_1v_1 = f_1.
\end{equation}
Since $(u,v)$ and $(p,q)$ are involutions, and the order of $(g,f)$ is three, 
we have the equalities 
\begin{equation}
\label{formula2}
u_0^2 = p_0^2 = v_1^2 = q_1^2 = 1, ~~~ \mbox{and} ~~~ g_0^3 = f_1^3 = 1.
\end{equation}
The required identities  $g_0 = f_1 = 1$ follow easily from 
(\ref{formula1}) and (\ref{formula2}).
\end{proof}

Next, take $(\alpha,\beta)\in  {\mathbb D} \times {\mathbb D}$ so that 
$$
l_m\bigl((\alpha,\beta)\bigl) = \Pi_m\bigl( (1 + \alpha t^m, t + \beta t^m)\bigr) 
\in \mu(K_4)\subseteq \mu(S_4) \subseteq {\R}_m
$$
has order 2. We claim that the product 
$$
\Pi_m\bigl( (1 + \alpha t^m, t + \beta t^m)\bigr) \Pi_m\bigl((g,f)\bigr)
$$
can not have order three in ${\R}_m$. Indeed, since 
$\Pi_m\bigl( (1 + \alpha t^m, t + \beta t^m)\bigr)$ has order 2 in ${\R}_m$, 
the equality 
$$
(1 + \alpha t^m, t + \beta t^m)^2 = (1 + 2\alpha t^m + h.o.t., ~ t + 2\beta t^m + h.o.t.)
$$
gives us $2\alpha = 2\beta = 0\in {\mathbb D}$. If 2 is not a zero divisor in 
$\mathbb D$, we can stop the proof here. Otherwise consider the product
$\bigl(1 + \alpha t^m,\, t + \beta t^m\bigr)\bigl(g,\,f\bigl)$. Since $f_1=1$ 
and $g_0 = 1$, we have $f(t+\beta t^m)=f(t)+\beta t^m + h.o.t.$, and
\begin{align*}
(1+\alpha t^m)g(t + \beta t^m)=& g(t+\beta t^m)+\alpha t^m g(t+\beta t^m)\\
=&g(t) +\beta g_1t^m+\alpha t^m + h.o.t.
\end{align*}
Therefore in ${\R}_m$ we can write 
$$
\Pi_m\Bigl((1 + \alpha t^m, t + \beta t^m)(g,\,f)\Bigr) = 
\Pi_m\Bigr(\bigl((1+\alpha t^m)g(t + \beta t^m),f(t+\beta t^m)\bigr)\Bigl)
$$
\begin{equation}
\label{order6}
= \bigl(g(t) + (\alpha + \beta g_1)t^m , ~ f(t) + \beta t^m \bigr)_m,
\end{equation}
where we used subindex $m$ in $\bigl(g,f\bigr)_m$ to denote the element 
$\Pi_m\bigl( (g,\, f )\bigr)$ in the truncated group ${\R}_m$. Next we take 
the element in (\ref{order6}) and raise it to the third power. We will work in the 
truncated group ${\R}_m$, and to keep the notations simple, we will not write 
the terms of the order higher than $m$ (and skip the ``$h.o.t.$"). Let us also 
denote the element in (\ref{order6}) by $C$. Then  
$$
C^2 = \bigl(g(t) + (\alpha + \beta g_1)t^m , ~ f(t) + \beta t^m \bigr)_m
\bigl(g(t) + (\alpha + \beta g_1)t^m , ~ f(t) + \beta t^m \bigr)_m
$$
$$
= \Bigl(\bigl(g(t) + (\alpha + \beta g_1)t^m\bigr)\cdot \bigl(g(f(t) + \beta t^m) + 
(\alpha + \beta g_1)(f(t) + \beta t^m)^m\bigr),~
$$
$$
f(f(t) + \beta t^m) + \beta (f(t) + \beta t^m)^m\Bigr)_m.
$$
Since $g_0 = f_1 = 1$, we have 
$\beta (f(t) + \beta t^m)^m = \beta f_1^mt^m = \beta t^m$ in ${\R}_m$, and using 
$2\beta = 0$ we obtain 
$$ 
f(f(t) + \beta t^m) + \beta (f(t) + \beta t^m)^m = f(f(t)) + \beta t^m + \beta t^m = f(f(t)).
$$
Similarly, in 
${\R}_m$, $g\bigr(f(t) + \beta t^m\bigl) = g(f(t)) + g_1\beta t^m$, so 
$$
g(f(t) + \beta t^m) + (\alpha + \beta g_1)(f(t) + \beta t^m)^m 
= g(f(t)) + \beta g_1 t^m + (\alpha +\beta g_1)t^m
$$
$$
= g(f(t)) + (\alpha  + 2 \beta g_1)t^m = g(f(t)) + \alpha t^m.
$$
Therefore (using $g_0 = 1$ with $2\alpha = 0$),
$$
\bigl(g(t) + (\alpha + \beta g_1)t^m\bigr)\cdot \bigl(g(f(t) + \beta t^m) + 
(\alpha + \beta g_1)(f(t) + \beta t^m)^m\bigr)
$$
$$
= \bigl(g(t) + (\alpha + \beta g_1)t^m\bigr)\cdot \bigl( g(f(t)) + \alpha t^m\bigr)
$$
$$
= g(t)g(f(t)) + \alpha t^m +  (\alpha + \beta g_1)t^m = g(t)g(f(t)) + \beta g_1 t^m,
$$
that is 
$$
C^2 = \Bigl(g(t)g(f(t)) + \beta g_1 t^m, f(f(t))\Bigr)_m.
$$
Furthermore, similar computations will show that 
$$
C^3 = \Bigl(g(t)g(f(t)) + \beta g_1 t^m, f(f(t))\Bigr)_m 
\Bigl(g(t) + (\alpha + \beta g_1)t^m , ~ f(t) + \beta t^m \Bigr)_m
$$
$$
= \Bigl( g(t)g(f(t))g(f(f(t))) + \alpha t^m, f(f(f(t))) + \beta t^m\Bigr)_m
$$
\begin{equation}
\label{lastone}
= \bigl(1 + \alpha t^m, t + \beta t^m\bigr)_m.
\end{equation}
Thus, in order $C^3 = \bigl(1, t\bigr)_m$ we must have $\alpha = \beta = 0$, which 
contradicts our choice of $(\alpha,\beta)$ producing $l_m\bigl((\alpha,\beta)\bigr)$ an 
element of order 2. This argument finishes the last possible case, and completes our  
proof of the theorem.
\end{proof}

We would like to note that impossibility of an embedding $S_n\to {\R}$ for $n\geq 5$ 
follows also from the corresponding statement for the alternating group of degree 
five $A_5$. The proof in this case is similar to the one above, but much shorter, 
due to simplicity of $A_5$. In fact, since the argument uses only the 
non-commutativity, finiteness, and simplicity of $A_5$, it is immediately generalized 
to all finite simple non-abelian groups, including many projective special linear 
groups over finite fields.

\begin{theorem}
Any finite non-abelian simple group $A$ cannot be embedded into the Riordan 
group ${\R}$ with coefficients in any commutative ring.
\end{theorem}
\begin{proof}
Assuming the contrary, there exists a monomorphism $\mu:A \to {\R}$. In this case, 
Proposition 2. guaranties existence of a monomorphism $\Pi_m\circ \mu: A \to {\R}_m$, 
where $m\geq 1$ is the smallest such natural number. Hence the composition 
$$
P_{m-1}\circ\Pi_m\circ \mu: A \longrightarrow {\R}_m \longrightarrow {\R}_{m-1}
$$
would have a nontrivial kernel, and since $A$ has no proper normal subgroups it
means that $\ker\bigl(P_{m-1}\circ\Pi_m\circ \mu \bigr) = A$. For $m\geq 2$ it
would imply $A\leq {\mathbb D}\times {\mathbb D}$, 
which is clearly impossible. For $m = 1$, we would get  
$A\leq {\mathbb D} \rtimes_{\varphi} {\mathbb D}^*$, 
which is also impossible using again the simplicity of the group $A$ and the fact that 
$ {\mathbb D} \rtimes_{\varphi} {\mathbb D}^*$ is a semidirect product of two 
abelian groups.
\end{proof}

\begin{corollary}
Let ${\mathbb F}$ be a finite field. If either $n > 2$, or $n = 2$ and $|{\mathbb F}| > 3$, 
then the projective special linear group ${\rm PSL}_n(\mathbb F)$ is not embeddable 
into ${\R}$.
\end{corollary}
\begin{proof}
It follows immediately from the result above and a theorem by Jordan
and Dickson, saying that under the given assumptions on $n$ and $\mathbb F$, 
the group ${\rm PSL}_n(\mathbb F)$ is simple 
(see \S 3.2.9 in Robinson \cite{Robinson}).
\end{proof}
Note that $A_5\cong {\rm PSL}_2({\mathbb F}_5)$, but in general the 
groups ${\rm PSL}_n(\mathbb F)$ are not of alternating type. For example 
${\rm PSL}_2({\mathbb F}_7)$, which is of order 168 and isomorphic to 
${\rm GL}_3({\mathbb F}_2)$ (see Brown and Loehr \cite{Brown}). Notice also 
that ${\rm PGL}_2(\inte)$, being a Coxeter group, is not embeddable 
into ${\R}(\comp)$ (see \S 2 in \cite{HeKrylov}).

Our Theorems 4 and 6 raise a question if the alternating group $A_4$, 
which has only one proper normal subgroup $K_4$, has a faithful representation by 
Riordan arrays. We end the discussion here with an example of such 
representation when the ground ring $\mathbb D$ has 2 as a zero 
divisor, and a primitive cube root of unity 
$\omega \in \mathbb D$. For instance, we 
can take $\mathbb D = \inte_6[\omega]$, where $\omega = (-1+\sqrt{-3})/2$. 
Here are more details. Consider the factor ring 
$$
\mathbb D = \inte_6[X]/\langle X^2 + X + 1\rangle,
$$
where $X$ is indeterminate over $\inte_6$ (the ring of integers modulo 6), and 
the principal ideal $\langle X^2 + X + 1\rangle$ is generated by the irreducible 
over $\inte_6$ polynomial $p(X) = X^2 + X + 1$. In particular, this ring $\mathbb D$ 
is finite, and all its elements can be represented by linear polynomials
$$
\mathbb D  \cong \left\{a + b\cdot X ~|~a,b\in\inte_6,~X^2 = 5 + 5 X\right\},
$$
with the multiplication given by 
$$
(a + b X)(p + qX) = (ap + 5bq) + (aq + bp + 5bq)X.
$$

In the following theorem we give an example of a particular
embedding of $A_4$ into the Lagrange subgroup ${\cal L}(\mathbb D)$, 
and also prove that $A_4$ can not be embedded into the substitution 
group ${\cal J}$ over any commutative ring.

\begin{theorem} 
A subgroup of the Lagrange subgroup 
${\cal L}(\mathbb D)$ generated by the Riordan arrays 
$$
u:=\bigl(1,\, t/(1 - 3t)\bigr), ~~ \mbox{and} ~~ 
w: = \bigl(1,\, X t\bigr)
$$
is isomorphic to $A_4$. Also, there exists no embedding of $A_4$ into the 
substitution group ${\cal J}$ with coefficients in an arbitrary commutative ring $R$.
\end{theorem}
\begin{proof}
Straightforward computations using 
$$
6 \equiv 1 + X + X^2 \equiv 0\pmod {\mathbb D},
$$ 
show that we have 
$$
u^2 = \bigl(1, \, t/(1 - 6t)\bigr) = (1, \, t), ~~~~~ w^3 = \bigl(1, \, X^3t\bigr) = \bigl(1, \, t\bigr),
$$
and  
$$
uw = \bigl(1, \, Xt/(1 - 3t)\bigr) , ~~~~~ (uw)^3 = (1,\, t).
$$
Furthermore, the following 12 
elements are all distinct in ${\L}_2(\mathbb D)$,
$$
u = \begin{pmatrix}
1 & 0 & 0\\
0 & 1 & 0\\
0 & 3 & 1\\
\end{pmatrix}, w = \begin{pmatrix}
1 & 0 & 0\\
0 & X & 0\\
0 & 0 & X^2\\
\end{pmatrix}, w^2 = \begin{pmatrix}
1 & 0 & 0\\
0 & X^2 & 0\\
0 & 0 & X\\
\end{pmatrix},
$$
$$
u^2 = w^3 = e = \begin{pmatrix}
1 & 0 & 0\\
0 & 1 & 0\\
0 & 0 & 1\\
\end{pmatrix}, uw = \begin{pmatrix}
1 & 0 & 0\\
0 & X & 0\\
0 & 3X & X^2\\
\end{pmatrix}, wu = \begin{pmatrix}
1 & 0 & 0\\
0 & X & 0\\
0 & 3X^2 & X^2\\
\end{pmatrix},
$$
$$
uw^2 = \begin{pmatrix}
1 & 0 & 0\\
0 & X^2 & 0\\
0 & 3X^2 & X\\
\end{pmatrix}, wuw = \begin{pmatrix}
1 & 0 & 0\\
0 & X^2 & 0\\
0 & 3 & X\\
\end{pmatrix},uwu = \begin{pmatrix}
1 & 0 & 0\\
0 & X & 0\\
0 & 3 & X^2\\
\end{pmatrix},
$$
$$
(uw)^2 = \begin{pmatrix}
1 & 0 & 0\\
0 & X^2 & 0\\
0 & 3(1+X^2) & X\\
\end{pmatrix}, 
wuw^2 = \begin{pmatrix}
1 & 0 & 0\\
0 & 1 & 0\\
0 & 3X & 1\\
\end{pmatrix}, w^2uw = \begin{pmatrix}
1 & 0 & 0\\
0 & 1 & 0\\
0 & 3X^2 & 1\\
\end{pmatrix}.
$$
Therefore, in terms of generators and relations, the subgroup generated by 
$u$ and $w$ has the presentation 
\begin{equation}
\label{PrA4}
\langle u, w ~|~ u^2 = w^3 = (uw)^3 = e\rangle,
\end{equation}
which is a group presentation of $A_4$ with the correspondences 
$(12)(34) \leftrightarrow u$, $(123) \leftrightarrow w$, and $(134) \leftrightarrow uw$.

To prove the second statement, let us assume to the contrary that there exists a 
monomorphism $\mu: A_4 \to {\cal J}$. In particular, using (\ref{PrA4}), 
it means that there exist two f.p.s. 
$$
f(t) = t + f_2t^2 + f_3t^3 +h.o.t. \in R[[t]] ,~~~ h(t) = t + h_2t^2 + h_3t^3 + h.o.t. \in R[[t]]
$$
such that 
\begin{equation}
\label{assume1}
\bigl(1, f(t)\bigr)^2 = \bigl(1, h(t)\bigr)^3 = \bigl(1, h(f(t))\bigr)^3 = (1,t).
\end{equation}
Let $f_r,r\geq 2$ be the smallest nonzero coefficient of $f(t)$, and 
$h_k,k\geq 2$ be the smallest nonzero coefficient of $h(t)$. In other words,
$$
f(t) = t + f_rt^r + h.o.t.~~~ \wedge ~~~ h(t) = t + h_kt^k + h.o.t., ~~ 
{\mbox s.t.} ~~ f_r\neq 0 \wedge h_k\neq 0.
$$
Then
$$
f(f(t)) = t + 2f_rt^r + h.o.t~~~ \wedge ~~~h(h(h(t))) = t + 3h_kt^k + h.o.t.,
$$
and (\ref{assume1}) imply $2f_r = 3h_k = 0$ in $R$. Furthermore, 
\begin{equation}
\label{HFt}
h(f(t)) = f(t) + h_kf^k(t) + h.o.t. = t + f_rt^r +\cdots + h_k(t + f_rt^r)^k + \cdots.
\end{equation}
If we suppose $r < k$, then (\ref{HFt}) together with (\ref{assume1}) imply $3f_r = 0$. 
Since $2f_r = 0 ~ \wedge ~3f_r = 0  \Longrightarrow f_r = 0$, which contradicts our 
assumption about $f_r$, we must have $k\leq r$. In this case we have
$$
[t^n]h(f(t)) = [t^n]h(t) = h_n ~\mbox{if}~ n < r, ~ \mbox{and} ~ [t^r]h(f(t)) = h_r + f_r.
$$
Let us denote the $n$-th degree Taylor polynomial of a f.p.s. $g(t)\in R[[t]]$ as 
$g_n(t)\in R[t]$. Then we can write 
$$
h(f(t)) = h_{r-1}(t) + (h_r + f_r)t^r + h.o.t = h_r(t) + f_rt^r + h.o.t.
$$ 
Writing only Taylor polynomials of degree $r$ for the 2nd and 3rd compositional powers of 
$h(f(t))$ we correspondingly obtain
$$
\Bigr(h(f(t))\Bigl)^{\circ 2}_r = h_r\bigl(h_r(t) + f_rt^r\bigr) + f_rh_r^r(t) = h^{\circ 2}_r(t) + 2f_rt^r,
$$
and 
$$
\Bigr(h(f(t))\Bigl)^{\circ 3}_r = h^{\circ 2}_r\bigl(h_r(t) + f_rt^r\bigr) + 2f_rh_r^r(t) = 
h^{\circ 3}_r(t) + 3f_rt^r,
$$
where $f^{\circ n}$ means applying the function $f$ to itself $n$ times. 
Since $h^{\circ 3}(t) = t$ and $\Bigr(h(f(t))\Bigl)^{\circ 3} = t$, we deduce that again 
$3f_r = 0$, and hence $f_r = 0$, which contradicts our assumption on $f_r$.
\end{proof}







\biboptions{longnamesfirst,square,semicolon}

\end{document}